\documentclass{amsart}
\usepackage{amsmath, amscd, amssymb, latexsym}

\newtheorem{theorem}{Theorem}[section]
\newtheorem{lemma}[theorem]{Lemma}

\newtheorem{proposition}[theorem]{Proposition}

\newtheorem{definition}[theorem]{Definition}
\newcommand{\mn}{\par\medskip\noindent}
\newcommand{\cB}{\mathcal B}
\newcommand{\cD}{\mathcal D}
\newcommand{\cF}{\mathcal F}
\newcommand{\cH}{\mathcal H}
\newcommand{\oH}{\overline{H}}
\DeclareMathOperator{\id}{id}
\def\Cz{\mathbb{C}}
\def\Nz{\mathbb{N}}
\newcommand{\hot}{\mathbin{\hat{\otimes}}}

\title{A remark about the Connes fusion tensor product}
\author{Andreas Thom}
\address{Andreas Thom, Mathematisches Institut der Universit\"at G\"ottingen,
Bunsenstr. 3-5, D-37073 G\"ottingen, Germany}
\email{thom@uni-math.gwdg.de}
\urladdr{http://www.uni-math.gwdg.de/thom}

\begin{document}

\begin{abstract} We analyze the algebraic structure of the Connes fusion tensor product (CFTP) in the case of bi-finite Hilbert modules over a
von Neumann algebra $M$. It turns out that all
complications in its definition disappear if one uses the closely related bi-modules of bounded vectors. We construct an equivalence
of monoidal categories with duality between a category of Hilbert bi-modules over $M$ with CFTP and some natural category of bi-modules over $M$ with
the usual relative algebraic tensor product.
\end{abstract}
\maketitle

\section{Introduction}

First of all, we want to stress that this article does not contain any substantially 
new results about the Connes fusion tensor product (CFTP). We do
not claim any originality, all calculation are straight-forward. It is the accumulation and presentation of the results that might
be useful to those, who want to learn about the CFTP. We hope that the presented results will have applications towards
a better understanding of the properties of algebraic quantum field theories, as constructed in \cite{Wa}. In particular, we hope
that the algebraization of the CFTP will give a suitable framework to develop a concrete model for elliptic cohomology, a process which
was started in \cite{TeSt}.

In this note we will explain some facts about the Connes fusion tensor product (CFTP)
of bi-finite Hilbert bi-modules over von Neumann algebras. We stick to the case of bi-modules which are left-modules and right-modules over the same algebra $M$. 
The more general case can be treated similarily. We denote the relative algebraic tensor product by $\otimes_M$. The CFTP will
be denoted by $\hot_M$ (see section \ref{CFTP} for a definition).

The organization of the article is as follows.
Section \ref{prel} contains preliminaries and definitions which are needed for the exposition. In Section \ref{dual} we recall the notion of conjugate and dualizable 
bi-module. Section \ref{theosec} contains Theorem \ref{maintheo}, which subsumes the main results. It states that the category of bi-finite dualizable Hilbert
bi-modules is equivalent to the category of bi-finitely generated projective bi-modules over $M$, respecting the monoidal structure and duality. Parts of
the proof are carried out in Section \ref{CFTP}, which also recalls the definition of the CFTP. Finally, Section 
\ref{conclusio} contains some concluding remarks.
\section{Preliminaries} \label{prel}
\subsection{Hilbert modules and Hilbert bi-modules}

Throughout the article, let $M$ be a von Neumann algebra with separable pre-dual. We denote by $M^{op}$ the von Neumann algebra, which is defined to be $M$ as complex vector space, but
with the opposite multiplication. A Hilbert space $H$, which is at the same time a left-module over $M$, is called a 
{\it Hilbert left-module} over $M$, if the module action is induced by a normal, $*$-preserving, and unital homomorphism $M \to \cB(H)$. 
A Hilbert right-module is defined similarily with $M^{op}$ in place of $M$.
A {\it Hilbert bi-module} over $M$ is a Hilbert space $H$ which is at the same time a Hilbert left-module and a Hilbert 
right-module over $M$, such that the actions of $M$ and
$M^{op}$ commute.
 
Let $\phi: M \to \Cz$ be a faithful and normal state. The assignment $(x,y) \mapsto \phi(y^*x)$ makes $M$ into a pre-Hilbert space, whose
completion we denote by $L^2_{\phi}(M)$. We denote the image of an element $x\in M$ under the natural inclusion $M \subset L^2_{\phi}(M)$ by
$\hat{x} \in L^2_{\phi}(M)$. Furthermore, we use the convention $\Omega = \hat 1$. 

The Hilbert space $L^2_{\phi}(M)$ is naturally a Hilbert left-module, since
$$\|y \rhd \hat{x}\|^2_{L^2_\phi(M)} = \|\hat{yx}\|^2_{L^2_\phi(M)} = \phi(x^*y^*yx) \leq \|y\|^2 \phi(x^*x) = \|y\|^2\|x\|^2_{L^2_\phi(M)},$$
so that left-multiplication by $y$ extends to a bounded operator on the completion. 

In Tomita-Takesaki theory (see \cite[Ch. VI-VIII]{Tk}), an operator $S$ is considered. It is defined as the closure of the conjugate linear assignment
$$\hat{b} \mapsto \hat{b^*}.$$ In general $S$ is an unbounded operator. Its polar decomposition $S = J \Delta^{1/2}$ is such that
$J$ is a conjugate-linear isometry and $\Delta^{1/2}$ is an unbounded, positive operator.

We define a right-module structure by setting $\xi \lhd b = J(b^*\rhd (J\xi))$. It follows from Tomita-Takesaki theory that the left
and right-module structure are commutants of each other, so that in particular, they give $L^2_{\phi}(M)$ the structure of a Hilbert bi-module.  

The Hilbert space $L^2_{\phi}(M)^{\oplus n}$ is at the same time a Hilbert left-module over $M$ and a Hilbert right-module over $M_n(M)$. Again, these actions
are commutants of each other. A similar remark applies to the obvious left-module structure over $M_n(M)$ and the diagonal action of $M$ on the
right. We will freely use these actions without further notational complication.

A Hilbert left-module $H$ is called {\it finitely generated}, if there is a bounded surjection
of left-modules from $L^2_{\phi}(M)^{\oplus n}$ to $H$, for some $n \in \Nz$. By polar decomposition, we can assume that the surjection is
a partial isometry. A similar definition applies to right-modules.
A Hilbert bi-module is called {\it bi-finite}, if it is finitely generated as left-module and as right-module.

By definition, all finitely generated Hilbert left-modules
over $M$ are unitarily equivalent to ${L^2_{\phi}(M)}^{\oplus n} \lhd p$, for some $n\in \Nz$ and some projection $p \in M_n(M)$.

The algebra $M$ is itself a bi-module over $M$, but usually not a Hilbert space. However, it is finitely generated projective as left and as right
module. Recall: a module is projective, if it is a direct summand in a free module.
In the sequel, we will encounter finitely generated projective
modules over $M$, which are associated to bi-finite Hilbert modules over $M$.

Let $H$ and $L$ be Hilbert left-modules (resp. right-modules). We denote by \begin{eqnarray*} 
\hom_{M-}(H,L) \quad \mbox{(resp.} \quad \hom_{-M}(H,L)) \end{eqnarray*} the vector space of
bounded left-module maps (resp. right module maps). 

Let $H$ and $L$ be either (Hilbert) left-modules, right-modules or bi-modules over $M$. 
The algebraic direct sum of the vectorspaces $H$ and $L$ carries obvious module structures, such that
$H \oplus L$ is categorical sum in the categories of (Hilbert) left-modules, right-modules and bi-modules over $M$, respectively.

\subsection{Modules of homomorphisms}

Let $H$ and $L$ be Hilbert bi-modules. We consider $\hom_{-M}(H,L)$ and $\hom_{M-}(H,L)$. There are natural 
bi-module structures on these vector spaces. For $x \in \hom_{-M}(H,L), b \in M$ and $\xi \in H$, we define 
\begin{eqnarray*} (b \rhd x)(\xi) = b \rhd x(\xi) \quad \mbox{and} \quad (x \lhd b)(\xi) = x(b \rhd \xi).
\end{eqnarray*}
Similarily, for $x \in \hom_{M-}(H,L), b \in M$ and $\xi \in H$, we define
\begin{eqnarray*}
(b \rhd x)(\xi) = x(\xi \lhd b) \quad \mbox{and} \quad (x \lhd b)(\xi) = x(\xi) \lhd b.
\end{eqnarray*}
It is an easy check, that these formula indeed define bi-module structures. 
The adjoint of $x \in \hom_{-M}(H,L)$ (resp. $x \in\hom_{M-}(H,L)$) is denoted by $x^* \in \hom_{-M}(L,H)$ (resp.
$x^* \in \hom_{M-}(L,H)$). It 
satisfies $(x(\xi),\eta)_L = (\xi,x^*(\eta))_H$. Concerning the module structure, we have (in either case)
\begin{eqnarray}\label{adjoint} (b \rhd x)^* = x^* \lhd b^* \quad  \mbox{and} \quad (x \lhd b)^* = b^* \rhd x^*. \end{eqnarray}
Note that there is a canonical isomorphism $$\alpha: \hom_{-M}(L^2_{\phi}(M),L^2_{\phi}(M)) \to M$$ of von Neumann algebras
and bi-modules over $M$. We have $\alpha(x) \rhd \xi = x(\xi)$, for $x \in \hom_{-M}(L^2_{\phi}(M),L^2_{\phi}(M))$. Similarily, there
is a canonical isomorphism $$\alpha': \hom_{M-}(L^2_{\phi}(M),L^2_{\phi}(M)) \to M,$$ of bi-modules 
which satisfies $x(\xi) = \xi \lhd \alpha'(x)$. 

\begin{definition}
Let $H$ be a Hilbert left-module. We define the bi-module of left-bounded vectors by $\cD(H,\phi)=hom_{-M}(L^2_{\phi}(M),H)$.
\end{definition}

There is a similar definition for {\it right-bounded vectors} in terms of left-module maps. In the case of dualizable modules, we will show in Theorem
\ref{left-right}, that the left-bounded and right-bounded vectors are naturally isomorphic as bi-modules over $M$.

\section{Duality} \label{dual}
\subsection{Conjugate bi-modules}

Let $H$ be a bi-module over $M$. We denote by $\oH$ the {\it conjugate bi-module}. It is defined to be $H$ as a real linear space, but
equipped with the conjugate complex structure. Denote the image of $\xi \in H$ (under the real linear identification) 
in $\oH$ by $\overline{\xi}$. We define a bi-module
structure over $M$ by $b \rhd \overline{\xi} = \overline{\xi \lhd b^*}$ and $\overline{\xi} \lhd b = \overline{b^* \rhd \xi}$.
If $H$ is a Hilbert bi-module, then $\oH$ is naturally a Hilbert bi-module, taking $(\overline{\xi},\overline{\eta})_{\oH}=\overline{(\xi, \eta)_H}$.
Note that $\oH$ is bi-finite if and only if $H$ is bi-finite.

\begin{proposition} \label{natiso}
There is a natural unitary isomorphism of bi-modules over $M$ $$\beta: L^2_{\phi}(M) \stackrel{\sim}{\to} \overline{L^2_{\phi}(M)}.$$
\end{proposition}
\begin{proof} The map $\beta$ is given 
by $$\hat m \mapsto \overline{\Omega \lhd m^*}= \overline{J(\hat m)}.$$ It is unitary since $J$ is a conjugate-linear isometry. \end{proof}

In the sequel we provide some isomorphisms between modules of homomorphisms.

\begin{proposition} \label{iso1} Let $H$ and $L$ be Hilbert bi-modules over $M$.
We have a natural isomorphism $$\hom_{-M}(K,H) \stackrel{\sim}{\to} \overline{\hom_{-M}(H,K)}$$ of bi-modules over $M$.
Similarily, there is an isomorphism 
$$\hom_{M-}(K,H) \stackrel{\sim}{\to} \overline{\hom_{M-}(H,K)}$$ of bi-modules over $M$.
\end{proposition}
\begin{proof} The map is given by the assignment $x \mapsto \overline{x^*}$. Using equations \ref{adjoint}, it is easily checked that this is map of
bi-modules by formula. \end{proof}

\begin{proposition} \label{iso2} There is 
a natural isomorphism $$\hom_{M-}(H,K) \stackrel{\sim}{\to} \overline{\hom_{-M}(\overline{H},\overline{K})},$$
of bi-modules over $M$. \end{proposition}
\begin{proof} The map is given by $x \mapsto \overline{x^\dagger}$, where $x^\dagger(\overline{\xi})= \overline{x(\xi)}$.

The next computations show, that it is indeed an isomorphism of bi-modules.
We have $$(b \rhd x)^{\dagger}(\overline{\xi}) = \overline{b \rhd x(\xi)} =  \overline{x(\xi)} \lhd b^*,$$ so that
$$\overline{(b \rhd x)^{\dagger}} = \overline{x^\dagger \lhd b^*} = b \rhd \overline{x^\dagger}.$$ Similarily, 
$$(x \lhd b)^{\dagger}(\overline{\xi}) = \overline{(x \lhd b)(\xi)} = \overline{x( b \rhd \xi)} 
= \overline{x (\overline{\overline{\xi} \lhd b^*})},$$ so that
$$\overline{(x \lhd b)^{\dagger}} = \overline{b^* \rhd x^{\dagger}} = \overline{x^\dagger} \lhd b.$$
\end{proof}

Putting Proposition \ref{iso1} and Proposition \ref{iso2} together, we get a natural isomorphism of bi-modules as follows.
\begin{eqnarray} \hom_{-M}(L,H) \stackrel{\sim}{\to} \hom_{M-}(\oH,\overline{L}). \end{eqnarray}  

A case of special interest is the case $L=L^2_{\phi}(M)$. Composing
with the inverse of the isomorphism $\beta: L^2_{\phi}(M) \stackrel{\sim}{\to} \overline{L^2_{\phi}(M)}$ yields an identification
of bi-modules over $M$ as follows: \begin{eqnarray} \label{eqn1} \hom_{-M}(L^2_{\phi}(M),H) \stackrel{\sim}{\to} \hom_{M-}(\oH,L^2_{\phi}(M)).\end{eqnarray}

Making everything explicit, the image of $x \in \hom_{-M}(L^2_{\phi}(M),H)$ is denoted by $x'$ and $x'(\overline{\xi}) = Jx^*(\xi)$.

\subsection{Dual modules}

In order to speak about duality, we need a monoidal category (for definitions see \cite[Ch. VII]{ML}) . 
The category of Hilbert bi-modules over $M$ is equipped with the Connes
fusion tensor product (CFTP) which makes into a monoidal category. We will give a definition of the CFTP and a proof of the preceding assertion in 
section \ref{CFTP} and freely use its existence in this section. Section
\ref{CFTP} is independent of the results obtained in this section.

Let $H$ be a Hilbert bi-module. A dual of $H$ is a Hilbert module $H'$, together with morphisms $\eta_H: H \hot_M H' \to L^2_{\phi}(M)$
and $\varepsilon_H: L^2_{\phi}(M) \to H' \hot_M H$, which satisfy the usual adjointness relations, i.e.
\[  (\id_{H}\hot_M \varepsilon) \circ (\eta_H \hot_M \id_H) = \id_H \quad \mbox{and} \quad
 (\varepsilon \hot_M \id_{H'}) \circ (\id_{H'} \hot_M \eta_H) = \id_{H'}.\]

We call a Hilbert bi-module $H$ {\it dualizable}, if $\oH$ is a dual of $H$. (A similar definition applies to bi-modules over $M$ with
respect to the relative algebraic tensor product.)

The notion of duality is self-dual. Thus, if $H$ is dualizable, then so is $\overline{H}$, we have $\varepsilon_{\oH}={\eta_H}^*$ 
and $\eta_{\oH}= {\varepsilon_H}^*$. 

The duality gives a natural identification of bi-modules 
\begin{eqnarray} \label{eqn2} \hom_{M-}(L^2_{\phi}(M),H) \stackrel{\sim}{\to}\hom_{M-}(\oH,L^2_{\phi}(M)) \end{eqnarray} 
via the assignment $\phi \mapsto \eta_H \circ (\phi \otimes \id_{H'}).$ This
yields the following theorem.

\begin{theorem} \label{left-right}
If $H$ is dualizable, then there is a natural isomorphism of bi-modules 
$$\hom_{-M}(L^2_{\phi}(M),H) \stackrel{\sim}{\to} \hom_{M-}(L^2_{\phi}(M),H).$$ 
\end{theorem}
\begin{proof} This follows from the isomorphisms in formula \ref{eqn1} and \ref{eqn2}. \end{proof}

\section{Main theorem} \label{theosec}
\subsection{Finiteness}

There is a natural linear inclusion map $\cD(H,\phi) \mapsto H$ which is given by $x \mapsto x(\Omega)$. It is a map of left modules over $M$, but
in general not a map of right modules. 

If $H$ is a Hilbert bi-module, finitely generated as a Hilbert right-module over $M$, it is isomorphic to $p \rhd L^2_{\phi}(M)^{\oplus n}$ as a right
module, for some $n \in \Nz$ and projection $p \in M_n(M)$. We easily get that $\cD(H,\phi) = pM^{\oplus n}$ as right-modules. Indeed, this
follows since $H \mapsto \hom_{-M}(L^2_{\phi}(M),H)$ (seen as a right-module) is functorial for right-module maps and preserves direct sums.
Thus, $\cD(H,\phi)$ is finitely generated projective as a right-module.

If $H$ is dualizable, we find a similar argument showing that $\cD(H,\phi)$ is finitely generated projective as a left-module. Indeed, this follows from
Theorem \ref{left-right}.
\begin{proposition} \label{bfgp}
If $H$ is a dualizable bi-finite Hilbert module over $M$, then $\cD(H,\phi)$ is finitely generated projective as left-module and as right-module over $M$.
\end{proposition}

\begin{definition}
A bi-module over $M$ which is finitely generated projective as a right-module and as a left-module is called a {\it bfgp bi-module} over $M$.
\end{definition}

The term 'bfgp' encodes the words {\it bi-finitely generated projective}.

\begin{definition}
We denote by $\cH_M$ the category of Hilbert bi-modules over $M$ and by $\cH^d_M$ the category of 
bi-finite dualizable Hilbert bi-modules over $M$. The morphisms in $\cH_M$ and $\cH^d_M$ are defined to be bounded maps of bi-modules.
\end{definition}

\subsection{Inner products}

Let $H$ be a bi-finite Hilbert bi-module over $M$. The bi-modules of the form $\cD(H,\phi)$ have an additional 
feature, they carry a $M$-valued inner product, which is given by the formula $$(x,y) = \alpha(y^*x).$$ In general, 
a $M$-valued inner product is a real linear map $H \times H \to M$, which is complex
linear in the first variable and satisfies $$(x,y) = (y,x)^* \quad \mbox{and} \quad (x,x) \geq 0, $$for all $x,y \in H$.

A $M$-valued inner product of a bi-module $H$ is called {\it complete} if the norm $x \mapsto \| (x,x)^{1/2} \|$ is complete. 
It is called {\it compatible} with the module actions, if 
\begin{eqnarray*} 
(b \rhd x,y) &=& (x,b^* \rhd y) \\
(x \lhd b,y) &=& (x,y) b.
\end{eqnarray*}

\begin{proposition} \label{exinn}
Let $M$ be a bfgp bi-module over $M$. There exists a complete and compatible $M$-valued inner product on $M$.
\end{proposition}
\begin{proof}
Let $H$ be a bfgp bi-module over $M$. Again, $H$ is isomorphic to $pM^{\oplus n}$ as right-module,
for some $n \in \Nz$ and $p \in M_n(M)$. We define an inner product on $H$ by setting 
$$((b_i)_{i=1}^n,({b'}_{i=1}^n))= \sum_{i=1}^n {{b'}_i}^*b_i.$$ 
One easily checks that this inner product is complete and compatible. However, a choice was involved. 
\end{proof}

\begin{proposition} \label{summa}
Let $H$ be a bi-finite dualizable Hilbert bi-module. The bi-module $\cD(H,\phi)$ is bfgp and the canonical $M$-valued 
inner product on $\cD(H,\phi)$ is complete and compatible.
\end{proposition}
\begin{proof} The first part is Proposition \ref{bfgp}. Compatibility and completeness are easily checked. \end{proof}

If $H$ and $K$ are bi-modules over $M$, which carry an inner product, then there is a natural inner product
on $H \otimes_M K$, given by the formula 
$$(x_1 \otimes y_1,x_2 \otimes y_2)_{H \otimes_M K} = ((x_1,x_2)_H \rhd y_1,y_2)_K.$$

\begin{proposition}
Every bfgp bi-module over $M$ is dualizable.
\end{proposition}
\begin{proof}
The duality is induced by the inner product, but let us first start with a plain bfgp bi-module $H$ over $M$ and let us show that
$\hom_{-M}(H,M)$ is a dual of $H$.
First of all, there is a map $\eta_H: \hom_{-M}(H,M) \otimes_M H \to M$, given by $\phi \otimes h \mapsto \phi(h)$.

The other structure map of the duality $\varepsilon_M :M \to H \otimes_M \hom_{-M}(H,M)$ is described by some central element in
$H \otimes \hom_{-M}(H,M)$. Note, that $H$ embeds in $M^{\oplus n}$ as right-module, with projection $p \in M_n(M)$, such
that
$$H \otimes_M \hom_{-M}(H,M) \subset M^{\oplus n} \otimes_M \hom_{-M}(M^{\oplus n},M) \stackrel{\sim}{\to} M^{\oplus n^2}$$
Here, $H \otimes_M \hom_{-M}(H,M)$ identifies, as a complex vector space, with $pM_n(M)p$ and the bi-module structure 
is given by {\it one} normal, unital and $*$-preserving homomorphism $\gamma: M \to pM_n(M)p$. Clearly, $p \in pM_n(M)p$ 
is a central element, so that a morphism of bi-modules $\varepsilon_M: M \mapsto H \otimes_M \hom_{-M}(H,M)$ can be
defined. One easily checks, that the adjointness conditions are fulfilled, so that $H$ is indeed a dual of $\hom_{-M}(H,M)$.

Consider now a choice of a $M$-valued inner product on $H$ which exists by Proposition \ref{exinn}. 
We get a natural map of bi-modules $\oH \mapsto \hom_{-M}(H,M)$,
which is given by $\overline{k} \mapsto \{H \ni h \mapsto (h,k)_H \in M\}$. One easily shows that it is an isomorphism. Thus, $H$ is
a dual of $\oH$, for any bfgp bimodule. This implies the assertion.
\end{proof}
\subsection{Categories of bi-modules}

\begin{definition}
The category of bi-modules over $M$ is denoted by $\cB_M$.
Denote the category of bfgp bi-modules over $M$ by $\cF_M$. 
Furthermore, denote the category of bfgp bi-modules with a choice of
a complete and compatible inner product by $\cF^{*}_M$. In either case, morphisms are bi-module maps.
\end{definition}

If follows from Proposition \ref{exinn}, that the forgetful functor $\cF^{*}_M \to \cF_M$ is an equivalence of categories , 
although there is no canonical inverse equivalence. 

The assignment $H \mapsto \cD(H,\phi)$ extends to a functor $\cD_\phi: \cH_M \to \cB_M$ from the category of Hilbert modules over $M$ to
the category of bi-modules over $M$. We have shown in Proposition \ref{summa} 
that this functor descents to a functor $\cD_{\phi}: \cH^d_M \to \cF^*_M$, from
the category of bi-finite Hilbert bi-modules over $M$ to the category of bfgp bi-modules 
with complete and compatible inner products over $M$.

Let $H$ be a bfgp bi-module over $M$ with complete and compatible inner product. The relative tensor product 
$H \otimes_M L^2_{\phi}(M)$ is naturally a bi-finite Hilbert bi-module.
Its inner product is computed by the equation 
$$(h \otimes \xi ,h' \otimes \eta)= ((h,h')_H \rhd \xi, \eta)_{L^2_{\phi}(M)}.$$ At the moment we leave the assertion open, that
this construction produces dualizable modules. It will become apparent after the definition of CFTP. 
We will assume the assertion for the rest of the section.

\begin{lemma}
There is a functor $L^2_\phi: \cF^*_M \to \cH^d_M$, which sends $H \mapsto H \otimes_M L^2_{\phi}(M)$. 
\end{lemma}

We continue by proving two propositions that show that the functors $L^2_\phi$ and $\cD_\phi$ are actually equivalences of categories.

\begin{proposition} Let $H$ be a bi-finite Hilbert bi-module over $M$. There is a natural unitary equivalence
$$\cD(H,\phi) \otimes_M L^2_{\phi}(M) \stackrel{\sim}{\to} H $$ of Hilbert bi-modules over $M$.
\end{proposition}
\begin{proof}
First of all, we define a map $\cD(H,\phi) \otimes_{\Cz} L^2_\phi(M) \to H$ by sending $$x \otimes \xi \mapsto x(\xi).$$ This
clearly falls to a map from the relative tensor product over $M$. It is a map of left-modules, since
$b \rhd (x \otimes \xi) = (b \rhd x) \otimes \xi \mapsto (b \rhd x)(\xi) = b \rhd (x(\xi))$.

Again, the map is a unitary equivalence if $H$ is isomorphic to $L^2_{\phi}(M)^{\oplus n}$, and hence if 
$H$ is any finitely generated Hilbert module.
The assertion about the right-module structure follows similarily. \end{proof}

\begin{proposition}
Let $H$ be a bfgp bi-module with complete and compatible inner product. We have a natural isomorphism
$$H \stackrel{\sim}{\to} \hom_{-M}(L^2_{\phi}(M),H \otimes_M L^2_{\phi}(M))$$ of
bfgp bi-modules with inner product.
\end{proposition}
\begin{proof}
Let $h$ be an element in $H$. We define a map $h': L^2_{\phi}(M) \to H \otimes_M L^2_{\phi}(M)$ by sending $\xi \mapsto h \otimes \xi$. The
assignment $h \mapsto h'$ is clearly a map of
right-modules. Again, it is an isomorphism if $H= M^{\oplus n}$, as right-modules, and hence for 
any bfgp bi-module by functoriality of the above assignment for right-modules.
We still have to give an argument that shows that the inner products match. First of all, we compute the adjoint of $h'$. We claim that
${h'}^*(k \otimes \eta)= (h,k)_M^* \rhd \eta$, for $k \in H$ and $\eta \in L^2_{\phi}(M)$.
\begin{eqnarray*}
(h'(\xi), k \otimes \eta)_{H \otimes_M L^2_{\phi}(M)} &=& (h \otimes \xi, k \otimes \eta)_{H \otimes_M L^2_{\phi}(M)} \\
&=& ((h,k)_H \rhd \xi, \eta)_{L^2_{\phi}(M)} \\
&=& (\xi, {(h,k)_H}^* \rhd \eta)_{L^2_{\phi}(M)} \\
&=& (\xi, h'(k \otimes \eta))_{L^2_{\phi}(M)}
\end{eqnarray*}

But clearly, \begin{eqnarray*} (h'(\xi),k \otimes \eta)_{H \otimes_M L^2_{\phi}(M)} &=&
(\xi, {h'}^*k(\eta))_{L^2_{\phi}(M)} \\
&=& (\xi, \alpha({k'}^*h')^* \rhd \eta)_{L^2_{\phi}(M)},\end{eqnarray*}
so that $\alpha({k'}^*h')= (h,k)_M$ and the assertion follows. \end{proof}

The above propositions imply the parts of the following theorem.

\begin{theorem} \label{maintheo} The functors $L^2_{\phi}: \cF^*_M \to \cH^d_M$ and $\cD_\phi: \cH^d_M \to \cF^*_M$ constitute an equivalence of monoidal categories with duality.
\end{theorem}

We claim, that these functors extends to an equivalence of tensor categories. The above propositions show that these functors are at least
equivalences of categories. In order to make the remaining assertion precise, we recall the 
definition and some facts  about the CFTP. The proof Theorem \ref{maintheo} will be finished in section \ref{CFTP}.

The above result should be compared with Theorem $6.24$ in \cite{L}. Our methods allow to prove Theorem $6.24$ in the context of not necessarily
finite von Neumann algebras. However, since we stick to bi-modules Theorem \ref{maintheo} does not imply Theorem $6.24$ in \cite{L} and thus 
cannot be seen in complete anology. Since we know of no interesting applications of left-modules (rather that bi-modules) over 
infinite von Neumann algebras, we leave this issue aside.

\section{The Connes fusion tensor product} \label{CFTP}

\subsection{Definition}

The CFTP makes makes the category of Hilbert bi-modules over $M$  into a monoidal category. 
It is defined as follows. Let $H$ and $K$ be bi-finite Hilbert bi-modules over $M$.
$\cD(H,\phi) \otimes_{\Cz} \cD(K,\phi)$ carries a pre-inner product, given by 
$$(x_1 \otimes y_1, x_2 \otimes y_2) = \phi(\alpha(y_2^*   (\alpha (x_2^*  x_1) \rhd y_1))).$$ The CFTP of $H$ and $K$
over $M$ is defined to be, as a Hilbert space, the completion of $\cD(H,\phi) \otimes_{\Cz} \cD(K,\phi)$ with respect to the pre-inner product above. 
It is denoted by $H \hot_M K$.

The following computation shows that the {\it usual} algebraic tensor relations lie in the null-space of the inner product.
\begin{eqnarray*}
((x_1 \lhd b) \otimes y_1,x_2 \otimes y_2) &=& \phi(\alpha(y_2^*   (\alpha(x_2^*   (x_1 \lhd b)) \rhd y_1))) \\
&=&  \phi(\alpha(y_2^*   (\alpha(x_2^*   x_1) \rhd (b \rhd y_1))) \\
&=& (x_1 \otimes (b \rhd y_1),x_2 \otimes y_2)
\end{eqnarray*}

This implies that the CFTP is actually the completion of $\cD(H,\phi) \otimes_M \cD(K,\phi)$. The module structure on
$H \hot_M K$ still has to be specified. The left-module structure is the one that is extended by the completion. However, the right-module
structure is the one coming from $K$, rather that $\cD(K,\phi)$. The following alternative description of $H \hot_M K$ will make this
apparent.

In the literature, there are several other definitions and formulas for the inner product which 
have caused some confusion, in particular to those who wanted to understand the CFTP from a purely algebraic point of view. 
We want to review the definition of Connes and Sauvogeot \cite{Co, Sa} and the definition of Wassermann \cite{Wa}.

\subsection{The definitions of Connes and Wassermann}

In \cite{Co,Sa}, A. Connes and later J.-L. Sauvageot considered an inner product on $\cD(H,\phi) \otimes K$. It is given by the formula
$$(x_1 \otimes \xi_1,x_2 \otimes \xi_2)'= (\alpha(x_2^*   x_1) \rhd \xi_1,\xi_2)_K.$$

On vectors $\xi_i$, which are of the form $y_i(\Omega)$, for some $y_i \in \cD(K,\phi)$, this reads as follows:

\begin{eqnarray*}
(x_1 \otimes y_1(\Omega), x_2 \otimes y_2(\Omega))' &=& ((\alpha(x_2^*   x_1) \rhd y_1)(\Omega) , y_2(\Omega))_K \\
&=& ((y_2^*   ((\alpha(x_2^*   x_1) \rhd y_1)(\Omega),\Omega)_{L^2_{\phi}(M)} \\
&=& \phi(\alpha(y_2^*   (\alpha(x_2^*   x_1) \rhd y_1))) \\
&=& (x_1 \otimes y_1, x_2 \otimes y_2).
\end{eqnarray*}

We get that the inclusion $\cD(K,\phi) \otimes_{\Cz} \cD(H,\phi) \subset \cD(K,\phi) \otimes_{\Cz} H$ is isometric (with respect to
the pre-inner products) and hence the completions are isometric.

In \cite{Wa} A. Wassermann defines the CFTP of Hilbert bi-modules $H$ and $K$ over $M$ 
as a completion of $\hom_{-M}(L^2_{\phi}(M),H) \otimes_{\Cz} \hom_{M-}(L^2_{\phi}(M),K)$ with respect to the inner product

$$(x_1 \otimes y_1, x_2 \otimes y_2)'' = (\alpha(x_2^*x_1) \rhd \Omega \lhd \alpha'(y_2^*y_1),\Omega)_{L^2_{\phi}(M)},$$
which he calls a {\it four-point formula}.

The following computation implies that the completion of the natural morphism $\hom_{-M}(L^2_{\phi}(M),H) \otimes_\Cz 
\hom_{M-}(L^2_{\phi}(M),K) \to \cD(H,\phi) \otimes_\Cz K$ is isometric with respect to the pre-inner products. Furthermore,
it has dense image. Therefore, the completion is an isometric isomorphism.

\begin{eqnarray*} 
(x_1 \otimes y_1(\Omega), x_2 \otimes y_2(\Omega))' &=& (\alpha(x_2^*x_1) \rhd y_1(\Omega),y_2(\Omega))_{L^2_\phi(M)} \\
&=& (y_2^* (\alpha(x_2^*x_1) \rhd y_1(\Omega)),\Omega)_{L^2_\phi(M)} \\
&=& (\alpha(x_2^*x_1) \rhd y_2^*y_1(\Omega),\Omega)_{L^2_\phi(M)} \\
&=& (x_1 \otimes y_1, x_2 \otimes y_2)''
\end{eqnarray*}

Note, that in either case $H \hot_M K$ is {\it not} defined as a completion
of the algebraic tensor product $H \otimes_{\Cz} K$.

\subsection{Some properties}

From the definition above, it is already clear that
\begin{eqnarray*}
L^2_{\phi}(M) \hot_M H &\cong& H \\ H \hot_M L^2_{\phi}(M) &\cong& H,
\end{eqnarray*} and
\begin{eqnarray*} 
(H \oplus K) \hot_M L &\cong& (H \hot_M L) \oplus (K \hot_M L) \\
H \hot_M (K \oplus L) &\cong& (H \hot_M K) \oplus (H \hot_M L)
\end{eqnarray*} for all Hilbert bi-modules $H,K$ and $L$ over $M$.

\begin{proposition} \label{cdtens} Let $K$ and $H$ be bi-finite Hilbert bi-modules.
There is a natural isomorphism \begin{eqnarray} \label{assign} \cD(K,\phi) \otimes_M \cD(H, \phi) \to
\hom_{-M}(L^2_{\phi}(M),K \hot_M H) = \cD(K \hot_M H,\phi) \end{eqnarray} of bi-modules over $M$.
\end{proposition}
\begin{proof} 
The map is given by completing the map $M \to \cD(K,\phi) \otimes_M \cD(H, \phi)$, which
is associated to an element in $x \otimes y \in \cD(K,\phi) \otimes_M \cD(H,\phi)$. It is described by
$$m \mapsto x \otimes (y \lhd m),$$ and extends to $L^2_{\phi}(M)$ since
\begin{eqnarray*}
\|x \otimes (y \lhd m) \|^2 &=&\phi(m^*\alpha(y^*  (\alpha(x^*   x) \rhd y))m) \\
&\leq& \| \alpha(y^*  (\alpha(x^* x) \rhd y) \| \, \phi(m^*m) \\
&=& C \, \|\hat{m}\|^2_{L^2_{\phi}(M)}.
\end{eqnarray*}

The arrow in equation \ref{assign} is an isomorphism, if $K$ is isomorphic, as a right-module, to $L^2_{\phi}(M)^{\oplus n}$, 
and hence if $K$ is any bi-finite Hilbert bi-module over $M$, by functoriality for maps of Hilbert right-modules in the left variable. 
\end{proof}

Proposition \ref{cdtens} implies the associativity of the CFTP for bi-finite Hilbert bi-modules. Indeed, we can transport the associativity isomorphisms from 
the category $\cF^*$. 
It is now obvious that the equivalence of categories between $\cH_M$ and $\cF^*_M$ is an equivalence of monoidal categories. Indeed,
we have shown the existence of natural isomorphisms $\cD_{\phi}(K) \otimes_M \cD_\phi(H) \stackrel{\sim}{\to} \cD_{\phi}(K \hot_M H)$,
compatible with the monoidal structure maps. It follows by abstract reasoning, that the functor $L^2_\phi$ also preserves the monoidal structure. The next
proposition makes this explicit.

\begin{proposition}\label{L2tens}
There is a natural isomorphism
$$ K \otimes_M H \otimes_M L^2_{\phi}(M) \stackrel{\sim}{\to} (K \otimes_M L^2_{\phi}(M)) \hot_M (H \otimes_M L^2_{\phi}(M)),$$
compatible with the monoidal structure maps.
\end{proposition}
\begin{proof}
The map is given by $k\otimes h \otimes \hat{m} \mapsto k \otimes (h \lhd m)$. Again, we can show that
this map extends to an isomorphism of Hilbert bi-modules.
\end{proof}

In the proof of Theorem \ref{maintheo}, we left open the assertion, that the functors $\cD_\phi$ and $L^2_\phi$ 
preserve duality. We still have to show that the image of the functor $L^2_\phi$ consists of dualizable modules.

Let $H$ be a Hilbert bi-module and $\oH$ its dual. There is are natural isomorphisms of bi-modules over $M$

\begin{eqnarray*}
\overline{\cD(H,\phi)} &\cong& \overline{\hom_{-M}(L^2_{\phi}(M),H)} \\
&\cong& \hom_{M-}(\overline{L^2_{\phi}(M)},\oH) \\
&\cong& \hom_{M-}(L^2_{\phi}(M),\oH) \\
&\cong& \hom_{-M}(L^2_{\phi}(M),\oH).
\end{eqnarray*}

Here we used, Proposition \ref{iso2}, Proposition \ref{natiso} and Theorem \ref{left-right}. 

This shows, that the functor $\cD_\phi$ and hence $L^2_{\phi}$ preserve conjugacy of modules, up to natural isomorphism. We showed in Proposition \ref{cdtens} that it preserves the monoidal structure up to natural isomorphism. This implies that $L^2_{\phi}$ actually maps dualizable modules to dualizable modules. Since all
bfgp bi-modules are dualizable, the desired assertion is proven.

\section{Conclusion} \label{conclusio}

One might draw the conlcusion that all monoidal categories 
appearing in the realm of von Neumann algebras and satisfying the usual finiteness assumptions, 
can be understood from a completely algebraic point of view. This point is surely well-known 
to the experts in the field. We still believe that it will be helpful to those, that
start working in conformal field theory, subfactors (especially of von Neumann factors of typ III) etc. 
\mn
There is no reason to be scared by
the CFTP. This was maybe part of the confusion related to it.

\end{document}